\documentclass{amsart}
\usepackage{amssymb}
\usepackage{amsfonts}
\usepackage{amsmath}
\usepackage{graphicx}
\usepackage[pagebackref]{hyperref}
\newtheorem{thm}{Theorem}[section]
\newtheorem{cor}[thm]{Corollary}
\newtheorem{lem}[thm]{Lemma}

\theoremstyle{definition}

\newtheorem{exmp}[thm]{Example}

\theoremstyle{remark}

\DeclareMathOperator{\w.dim}{w.dim}

\newcommand{\field}[1]{\mathbb{#1}}

\newcommand{\N}{\field{N}}
\def\1{{\rm (1)}}
\def\2{{\rm (2)}}
\def\3{{\rm (3)}}
\def\4{{\rm (4)}}
\def\5{{\rm (5)}}

\begin{document}

\title[On Pr\"ufer-like conditions]{On Pr\"ufer-like conditions}

\author{C. Bakkari}

\address{Department of Mathematics, Faculty of Sciences and Technology, P. O. Box 2202,
University S. M. Ben Abdellah, Fez 30000, Morocco} \email{cbakkari@hotmail.com}

\date{July 10 2008}




\begin{abstract}
This paper deals with five extensions of the Pr\"ufer domain concept to commutative rings with zero divisors. We investigate the stability of these Pr\"ufer-like conditions under localization and homomorphic image. Our results generate new and original examples of Pr\"ufer-like rings.
\end{abstract}

\footnotetext[0]{ *Journal of Commutative Algebra, To appear
}\value{counter}\setcounter{footnote}{0}
\maketitle

\begin{section}{Introduction}

Throughout this paper all rings are commutative with identity
element and all modules are unital. In his article [22], Pr\"ufer introduced a new class of integral domains; namely
those domains in which all finitely generated ideals are invertible. Through the years, Pr\"ufer domains acquired a great many equivalent characterizations, each of which was extended to rings with zero-divisors in different ways. More precisely, we consider the following Pr\"ufer-like properties on a commutative ring $R$ [3, 4]:\\
\1 $R$ is semihereditary, i.e., every finitely generated ideal of $R$ is projective. \\
\2 The weak global dimension of $R$ is at most one. \\
\3 $R$ is arithmetical, i.e., every finitely generated ideal of $R$ is locally principal. \\
\4 $R$ is Gaussian, i.e., $c(fg) =c(f)c(g)$ for any polynomials $f, g\in R[X]$, where $c(f)$ is the content of $f$, that is, the ideal of $R$ generated by the coefficients of $f$.\\
\5 $R$ is Pr\"ufer, i.e., every finitely generated regular ideal is invertible (equivalently, every two generated regular ideal is invertible).

In [11], it is proved that each one of the above conditions implies the following
next one: (1) $\Rightarrow$ (2) $\Rightarrow$ (3) $\Rightarrow$ (4) $\Rightarrow$ (5). Also examples are given to show that,
in general, the implications cannot be reversed. Moreover, an investigation is carried out to see which conditions may be added to any of these properties in order to reverse the implications. Recall that in the domain context, the above five classes of Pr\"ufer-like
rings collapse to the notion of Pr\"ufer domain. From Bazzoni and Glaz [4, Theorem 3.12], we note that a Pr\"ufer ring $R$ satisfies anyone of the five conditions if and only if its total ring of quotients $Tot(R)$ satisfies the same condition. For more details on these notions, we refer the reader to [3, 4, 7, 10, 11, 13, 17, 21, 25].

This paper investigates the stability of the above five Pr\"ufer-like conditions under localization and homomorphic
image. In Section 2, we prove that if $R$ is a Pr\"ufer ring then so is $S^{-1}R$ for any multiplicative subset $S \subseteq R\setminus Z(R)$, where $Z(R)$ is the set of zero-divisors of $R$ (Theorem 2.1). We further show that the Pr\"ufer property is not, in general, stable under localization (Example 2.3). However, Theorem 2.5 asserts that if $R$ is Gaussian (resp., arithmetical, $\w.dim(R) \leq 1$, or semihereditary) then so is $S^{-1}R$ for any multiplicative subset $S$ of $R$.

In Section 3, Theorem 3.1 states that the homomorphic image of an arithmetical (resp., Gaussian) ring is arithmetical (resp., Gaussian). We also show that the remaining three Pr\"ufer-like conditions are not stable under homomorphic image (Examples 3.2 and 3.3). The section closes with a result (Theorem 3.4) which investigates the transfer of the five Pr\"ufer-like conditions to a particular case of homomorphic image; namely, direct product of rings.
\end{section}

\begin{section}{Localization of Pr\"ufer conditions}

In this section we present a detailed treatment of the
localization of the pre-mentioned Pr\"ufer-like conditions. We prove
that, unlike Pr\"ufer rings, the classes of Gaussian rings,
arithmetical rings, rings with $w.dim(R)\leq 1$ and semihereditary
rings are stable under localization.  We start with the following
theorem which states a condition under which the class of Pr\"ufer
rings is stable under localization.

\begin{thm}\label{locPr}
Let $R$ be a Pr\"ufer ring and $S$ a multiplicative subset of $R$ which is contained in $R\setminus Z(R)$.
Then $S^{-1}R$ is a Pr\"ufer ring.
\end{thm}

The proof will use the following Lemma.

\begin{lem}\label{LocPrlem}  Let $R$ be a Pr\"ufer ring and $S$ a multiplicative subset
of $R$ which is contained in $R\setminus Z(R)$.
Then $Tot(R) =Tot(S^{-1}R)$.
\end{lem}

\begin{proof} We claim that $R \subseteq S^{-1}R$. Indeed, let $0\not= a
\in R$ such that ${a\over{1}} ={0\over{1}}$. Then, there exists $t
\in S$ such that $ta =0$ and hence $a =0$ since $t \in S \subseteq
R\setminus Z(R)$. This means that $R \subseteq S^{-1}R$. Hence, every element
${a\over b}$ of $Tot(R)$, where $a \in R$ and $b \in R\setminus Z(R)$, can be
written as ${a\over b} ={{a/1}\over {{b/1}}}$ with
${b\over1} \in S^{-1}R \setminus Z(S^{-1}R)$. Thus, $Tot(R) \subseteq Tot(S^{-1}R)$. \\
Conversely, let $x ={{a/s}\over{a'/{s'}}} \in
Tot(S^{-1}R)$, where ${a\over s} \in S^{-1}R$ and ${a'\over s'} \in
S^{-1}R \setminus Z(S^{-1}R)$. We claim that $a's \in R\setminus Z(R)$. Indeed, let $b
\in R$ such that $ba's =0$. Then, ${({bs \over 1})({a'\over s'})}
={0\over 1}$ ($\in S^{-1}R$) and so ${{bs\over 1} ={0\over 1}}$
since ${a'\over s'} \in S^{-1}R \setminus Z(S^{-1}R)$. Thus, there exists $t
\in S$ such that $tsb =0$ and so $b =0$ since $ts \in S \subseteq
R\setminus Z(R)$. Hence, $x$ can be written as $x ={{a/1}\over {a'/1}}{{1/s}\over{1/s'}} = {{{(a/1)}{(1/s)}{(ss'/1}}\over {{(a'/1)}{(1/s')}{(ss'/1)}}} ={{(as')/1}\over
{(a's/1)}}$ since ${ss'\over 1} \in S^{-1}R \setminus Z(S^{-1}R)$.
Therefore, $x ={{as'\over a's}}$ since $a's \in R\setminus Z(R)$ which means
that $x \in Tot(R)$.
\end{proof}

\begin{proof}[Proof of Theorem~\ref{locPr}] 
One of the many characterizations of Pr\"ufer rings is that each overrings is integrally
closed. It is also clear by Lemma 2.2 that $R \subseteq S^{-1}R \subseteq Tot(R) (=Tot(S^{-1}R))$
for each subset $S \subseteq R \setminus Z(R)$. Thus $R$ Pr\"ufer implies $S^{-1}R$ is Pr\"ufer.
\end{proof}

The next example shows that the condition $S \subseteq R\setminus Z(R)$
cannot be dropped in Theorem \ref{locPr}. For this, we appeal to the
notion of trivial ring extension. We recall that for a ring $A$ and
an $A$-module $E$, the trivial ring extension of $A$ by $E$ (also
called the idealization of $E$ over $A$) is the ring $R:=
A\propto~E$ whose underlying group is $A \times E$ with
multiplication given by $(a_{1}, e_{1})(a_{2}, e_{2}) = (a_{1}a_{2},
a_{1}e_{2}+a_{2}e_{1})$. Considerable work, part of it summarized in
Glaz's book [9] and Huckaba's book [16], has been
concerned with trivial ring extensions. These have proven to be
useful in solving many open problems and conjectures for various
contexts in (commutative and non-commutative) ring theory, see for instance [9, 16, 18].

\begin{exmp}\label{Exm1}
Let $A =K[[X_{1},X_{2},X_{3}]] =K+M$ be a  power series ring over a
field $K$ and $M :=(X_{1}, X_{2}, X_{3})$. Let $E$
be an $A$-module such that $ME =0$ and let $R :=A \propto E$ be the
trivial ring extension of $A$ by $E$. Let $S$ be the multiplicative
subset of $R$ given by $S :=\{(X_{1},0)^{n}\ /\ n \in \N\}$ and $S_{0}$ the multiplicative subset of
$A$ given by $S_{0} :=\{X_{1}^{n}\ /\ n \in \N\}$. Then:
\begin{enumerate}
    \item $R$ is a Pr\"ufer ring.
    \item $S_{0}^{-1}A$ is a domain which is not Pr\"ufer.
    \item $S^{-1}R$ and $S_{0}^{-1}A$ are isomorphic rings. In particular, $S^{-1}R$ is not Pr\"ufer.
\end{enumerate}
\end{exmp}

 \begin{proof} \1 One may easily verify that $R$ is local with maximal ideal
$M \propto E$ and each element of $R$ is either a unit or a zero divisor.
Thus, $R =Tot(R)$ is Pr\"ufer.

\2 We clearly have $$S_{0}^{-1}A =S_{0}^{-1}K[[X_{1},X_{2},X_{3}]] =(S_{0}^{-1}K[[X_{1}]])[[X_{2},X_{3}]]
\subseteq qf(K[[X_{1}]])[[X_{2},X_{3}]].$$ Since $S_{0}^{-1}K[[X_{1}]]$ is Noetherian, by [9, Theorem 8.1.1]
$S_{0}^{-1}A$ is a domain with
$$\w.dim(S_{0}^{-1}A) = \w.dim(S_{0}^{-1}K[[X_{1}]])[[X_{2},X_{3}]] =\w.dim(S_{0}^{-1}K[[X_{1}]]) + 2 \geq 2.$$
In particular, $S_{0}^{-1}A$ is not a Pr\"ufer domain.

\3 Since $X_{1}E \subseteq ME =0$ and $X_{1} \in S_{0}^{-1}$,
then $S_{0}^{-1}E =0$. Thus, $S^{-1}(0 \propto E) =0$ and so
$S^{-1}R =\{{(a,0)\over (s,0)} \hspace{0.2cm} / \hspace{0.2cm}
a \in A$ and $s \in S_{0}\}$. Now, we easily check that:
$$f : S_{0}^{-1}A \longrightarrow  S^{-1}R $$
$${a\over s} \longmapsto  {{(a,0)\over (s,0)}}$$
is a ring isomorphism. In particular, $R$ is not a Pr\"ufer domain by \2.
\end{proof}

\begin{cor}\label{LocPr}
Let $R$ be a ring and $S$ a multiplicative subset of $R$ which is contained in $R\setminus Z(R)$. Then:
\begin{enumerate}
    \item If $R$ is Gaussian, then so is $S^{-1}R$.
    \item If $R$ is arithmetical, then so is $S^{-1}R$.
    \item If $\w.dim(R) \leq 1$, then $\w.dim(S^{-1}R) \leq 1$.
    \item If $R$ is semihereditary, then so is $S^{-1}R$.
\end{enumerate}
\end{cor}

\begin{proof}If $R$ satisfies one of the five Pr\"ufer-like conditions,
then so is $Tot(R) (=Tot(S^{-1}R))$ by [4, Theorem 3.12]. Also, $R$
is, in all cases, a Pr\"ufer ring and so $S^{-1}R$ is a Pr\"ufer
ring by Theorem 2.1. Therefore, $S^{-1}R$ satisfies the same
Pr\"ufer-like condition by [4, Theorem 3.12].
\end{proof}


The localization of a Pr\"ufer ring is not always a Pr\"ufer
ring by Example 2.3. For the other Pr\"ufer-like conditions, we
have:

\begin{thm}\label{locPrcond}
 Let $R$ be a ring and $S$ a multiplicative subset of $R$. Then:
\begin{enumerate}
    \item If $R$ is Gaussian, then so is $S^{-1}R$.
    \item If $R$ is arithmetical, then so is $S^{-1}R$.
    \item If $\w.dim(R) \leq 1$, then $\w.dim(S^{-1}R) \leq 1$.
    \item If $R$ is semihereditary, then so is $S^{-1}R$.
\end{enumerate}
\end{thm}

\begin{proof}  \1 Assume that $R$ is a Gaussian ring. Our aim is to
prove that for all polynomials $S^{-1}f$ and $S^{-1}g$, we have
$c_{S^{-1}R}(S^{-1}fS^{-1}g)
=c_{S^{-1}R}(S^{-1}f)c_{S^{-1}R}(S^{-1}g)$. Without loss of
generality, we may assume that $f
=\displaystyle\sum_{i=0}^{n}a_{i}X_{i}$ and $g
=\displaystyle\sum_{i=0}^{m}b_{i}X_{i} \in R[X]$. But, $c_{R}(fg)
=c_{R}(f)c_{R}(g)$ since $R$ is Gaussian. Hence:
\begin{eqnarray}
 \nonumber  c_{S^{-1}R}(S^{-1}f)c_{S^{-1}R}(S^{-1}g)&=& S^{-1}(c_{R}(f))S^{-1}(c_{R}(g)) \\
 \nonumber          &=& S^{-1}[c_{R}(f)c_{R}(g)] \\
\nonumber         &=& S^{-1}(c_{R}(fg)) \ (since\ R\ is\ Gaussian) \\
 \nonumber        &=& c_{S^{-1}R}(S^{-1}fS^{-1}g), {\rm as\hspace{0.1cm} desired.}
\end{eqnarray}

\2 Let $J$ be a finitely generated ideal of $S^{-1}R$ and $M$
be a maximal ideal of $S^{-1}R$. There exist a finitely generated
ideal $I$ of $R$ and a prime ideal $m$ of $R$ such that $J
=S^{-1}I$ and $M =S^{-1}m$. Next note that $S^{-1}R_{M}$ is naturally isomorphic
to $R_m$. As $R$ is arithmetical, $IR_{P}$ is principal for each prime $P$. It
follows that $J_{M} \cong IR_{m}$ is locally principal.

\3 and \4 are clear, completing the proof of Theorem 2.5.
\end{proof}


Now, we are able to construct a non-Gaussian Pr\"ufer ring.

\begin{exmp}\label{Exm2}
Let $R$ and $S$ be as in Example 2.3. Then:
\begin{enumerate}
    \item $R$ is a Pr\"ufer ring.
    \item $R$ is not a Gaussian ring.
     \end{enumerate}
\end{exmp}

 \begin{proof} \1 $R$ is a Pr\"ufer ring by Example 2.3(1).

\2 We claim that $R$ is not Gaussian. Deny, $S^{-1}R$ is a Gaussian ring by
Theorem 2.5(1). Then $S^{-1}R$ is a Pr\"ufer domain by Example 2.3,
which contradicts Example 2.3(3). Hence, $R$ is not a Gaussian ring.
\end{proof}

\end{section}

\begin{section}{Homomorphic image of Pr\"ufer conditions}

This  section studies the homomorphic image of Pr\"ufer-like rings. We show that the homomorphic image of a Gaussian (resp.,
arithmetical) ring is Gaussian (resp., arithmetical). On the other hand, we show that the other three classes of Pr\"ufer-like rings are not stable under homomorphic image (Examples 3.2 and 3.3).

\begin{thm}\label{homimpr}
\begin{enumerate}
    \item The homomorphic image of a Gaussian ring is Gaussian.
    \item The homomorphic image of an arithmetical ring is arithmetical.
     \end{enumerate}
\end{thm}

\begin{proof}
\1 By [12]. \\
\2 By [5].

\end{proof}

There have many examples presented of Pr\"ufer rings such that some homomorphic image
is not Pr\"ufer. One of the earliest sources for $R$ Pr\"ufer not implying $R/I$ Pr\"ufer
is the paper [6]. There are several such examples in the Example section of Huckaba's book [16]. \\
Now, we give a new example showing that the homomorphic image of a
Pr\"ufer ring is not, in general, a Pr\"ufer ring. This is also an example of
 a non-Gaussian Pr\"ufer ring.

\begin{exmp}\label{Exm3}
 Let $A$ be a non-valuation local domain, $M$ its maximal ideal, and $E$ an
 $A$-module with $ME =0$. Let $R :=A \propto E$ be the trivial ring extension of $A$ by $E$.
Consider the following ring homomorphism:
$$h :R   \longrightarrow  A$$
$$(a,e) \longmapsto  a.$$
Then:
\begin{enumerate}
    \item $R$ is a Pr\"ufer ring (since it is a total ring).
    \item $R$ is not a Gaussian ring by Theorem 3.1(1) since $h(R)=A$ is not a Gaussian domain.
     \end{enumerate}
\end{exmp}


 Now, we construct an arithmetical ring $R$ with $\w.dim(R) =\infty $ which shows that
 the homomorphic image of a semihereditary ring (resp., a ring with $\w.dim\leq 1$)
 is not necessarily semihereditary (resp., of weak dimension less than or equal to one).

\begin{exmp}\label{Exm4}
 Let $A =K[[X]]$ be the power series ring over a field $K$, $I =(X^{n})$ the ideal
of $A$ generated by $X^n$, where $n \geq 2$. Set $R :=A/I$. Then:
\begin{enumerate}
    \item $A$ is a discrete valuation domain.
    \item $R$ is an arithmetical Noetherian ring.
    \item $\w.dim(R) =\infty $. In particular, $R$ is not semihereditary.
     \end{enumerate}
\end{exmp}

 \begin{proof} \1 Trivial.

\2 $R$ is an arithmetical ring by Theorem 3.1(2) since $A$ is a Pr\"ufer domain. Also,
$R$ is Noetherian since the homomorphic image of Noetherian ring is Noetherian.

\3 Let $x^{i}$ be the image of $X^{i}$ in $R =A/I$. We denote by $(x^{i})$
the principal ideal of $R$ generated by $x^{i}$. It is easy to check that the following sequences:
$$0\longrightarrow  (x^{n-1})\longrightarrow  R \buildrel u \over\longrightarrow  (x)\longrightarrow  0$$
$$0\longrightarrow  (x)\longrightarrow  R \buildrel v \over\longrightarrow  (x^{n-1})\longrightarrow  0$$
where $u(r) =rx$ and $v(r) =rx^{n-1}$ for each $r \in R$,
are exact. But, the principal ideal $(x)$ of $R$ is not a projective
ideal since $xx^{n-1} =0$ and $R$ is local. Therefore, by the above
two exact sequences of $R$, we have $pd_{R}((x)) (=fd_{R}((x))
=\infty $ (since $R$ is Noetherian) which means that $\w.dim(R) =
\infty $.  In particular, $R$ is not semihereditary.
\end{proof}


Finally, we study a particular case of homomorphic images, that is,
the direct product of Pr\"ufer-like rings.

\begin{thm}\label{dirpr} Let $(R_{i})_{1\leq i\leq n}$ be a family of rings. Then:
\begin{enumerate}
\item $\displaystyle\prod_{i=1}^{n}R_{i}$ is Pr\"ufer if and only if so is
$R_i$ for each $i =1,\ldots, n$.
\item $\displaystyle\prod_{i=1}^{n}R_{i}$ is Gaussian if and only if so is
$R_i$  for each $i =1,\ldots, n$.
\item $\displaystyle\prod_{i=1}^{n}R_{i}$ is arithmetical if and only if so is
$R_i$ for each $i =1,\ldots, n$.
\item $\w.dim(\displaystyle\prod_{i=1}^{n}R_{i}) =\sup\{\w.dim(R_{i})\ |\ i =1, \ldots, n\}$.
\item $\displaystyle\prod_{i=1}^{n}R_{i}$ is semihereditary if and only if so is
$R_i$  for each $i =1,\ldots, n$.
\end{enumerate}
\end{thm}

\begin{proof} The proofs of the first three assertions are done by induction on $n$ and it suffices to check it for $n=2$.

\1 By [8].\\

\2 If $R_{1} \times R_{2}$ is a Gaussian ring, then, for each $i =1, 2$, $R_{i}$ is a Gaussian ring
  as homomorphic image of a Gaussian ring by Theorem 3.1(1).

Conversely, assume that $R_1$ and $R_2$ are Gaussian rings. Let $f
=\sum_{i=0}^{n}(a_{i},b_{i})X_{i}$ and $g
=\sum_{i=0}^{m}(c_{i},d_{i})X_{i}$ be two polynomials
in $(R_{1} \times R_{2})[X]$ and let $f_{1}
=\sum_{i=0}^{n}a_{i}X_{i} \in R_{1}[X]$, $f_{2}
=\sum_{i=0}^{n}b_{i}X_{i} \in R_{2}[X]$, $g_{1}
=\sum_{i=0}^{m}c_{i}X_{i} \in R_{1}[X]$, and $g_{2}
=\sum_{i=0}^{m}d_{i}X_{i} \in R_{2}[X]$. We will prove
that $c_{R_{1} \times R_{2}}(fg) =c_{R_{1} \times R_{2}}(f)c_{R_{1}
\times R_{2}}(g)$. First we note that $c_{R_{1} \times R_{2}}(f)
=(c_{R_{1}}(f_{1}),c_{R_{2}}(f_{2}))$. Hence, it is easy to see that
$c_{R_{1} \times R_{2}}(fg) =$ \\
$(c_{R_{1}}(f_{1}g_{1}),c_{R_{2}}(f_{2}g_{2}))
=(c_{R_{1}}(f_{1})c_{R_{1}}(g_{1}),c_{R_{2}}(f_{2})c_{R_{2}}(g_{2}))$
(since $R_1$ and $R_2$ are Gaussian rings) $=
(c_{R_{1}}(f_{1}),c_{R_{2}}(f_{2}))(c_{R_{1}}(g_{1}),c_{R_{2}}(g_{2})) =c_{R_{1} \times R_{2}}(f)c_{R_{1} \times R_{2}}(g)$.

\3 If $R_{1} \times R_{2}$ is an arithmetical ring, then, for each $i =1, 2$, $R_{i}$ is an arithmetical ring
 as homomorphic image of an arithmetical ring by Theorem 3.1(2).

Conversely, assume that $R_1$ and $R_2$ are arithmetical rings. Let
$J$ be a finitely generated ideal of the ring $R_{1} \times R_{2}$,
and $M$ a maximal ideal of $R_{1} \times R_{2}$. Write $J =I_{1}
\times I_{2}$, where $I_i$ is a finitely generated ideal of $R_i$,
and $M$ is either $m_{1} \times R_{2}$ or $R_{1} \times m_{2}$,
where $m_{i} \in Max(R_{i})$ for $i =1, 2$.
We may assume that $M =m_{1} \times R_{2}$ (the case $M =R_{1}
\times m_{2}$ is similar). We will prove that $J_{M} =(I_{1} \times
I_{2})_{m_{1} \times R_{2}}$ is principal. For this, consider an
element ${(a_{1},a_{2})\over(s_{1},s_{2})} \in J_{M}$, where $a_{i}
\in I_{i}$ for $i =1,2$ and $(s_{1},s_{2}) \in (R_{1} \times
R_{2})-(m_{1} \times R_{2})$. So, we have ${{(a_{1},a_{2})\over
(s_{1},s_{2})}} ={{(a_{1},a_{2})(1,0)}\over {(s_{1},s_{2})(1,0)}}
={(a_{1},0)\over (s_{1},0)}$ since $(1,0) \in (R_{1} \times
R_{2})-(m_{1} \times R_{2})$. Assume that $I_{1_{m_{1}}}
=aR_{1_{m_{1}}}$, where $a \in I_{1}$ since $R_1$ is supposed to be
arithmetical. Therefore, ${(a_{1},a_{2})\over (s_{1},s_{2})}
={(a_{1},0)\over (s_{1},0)} \in (a,0)(R_{1} \times R_{2})_{m_{1}
\times R_{2}}$. On the other hand, $(a,0) \in J$ imply that
$(a,0)(R_{1} \times R_{2})_{m_{1} \times R_{2}} \subseteq J_{M}
=J(R_{1} \times R_{2})_{m_{1} \times R_{2}}$. Finally, $J_{M}
=(a,0)(R_{1} \times R_{2})_{M}$ is principal and so
$R_{1} \times R_{2}$ is an arithmetical ring.

\4 and \5 are clear since a ring $R$ is semihereditary
if and only if it is coherent and $\w.dim(R) \leq 1$, a finite direct
product of rings is coherent if and only if each component is
coherent, and $\w.dim(\prod_{i=1}^{n}R_{i})
=sup\{\w.dim(R_{i})\ |\ i =1, \ldots, n\}$.
\end{proof}


We close this paper by constructing a non-total non-Gaussian
Pr\"ufer ring.

\begin{exmp}\label{Exm5}
Let $R$ be as in Example 2.3, $T$ a Pr\"ufer domain which is not a field, and
$L :=R \times T$ the direct product of $R$ by $T$. Then:
\begin{enumerate}
\item $L$ is a non-total ring since $(1,a) (\in L)$ is neither unit nor zero-divisor in $L$
    for each non invertible element $a \in T$.
\item $L$ is a Pr\"ufer ring by Theorem 3.4(1) since $R$ and $T$ are Pr\"ufer rings.
\item $L$ is not a Gaussian ring by Theorem 3.4(2) since $R$ is not a Gaussian ring.
\end{enumerate}
\end{exmp}
\end{section}



\begin{thebibliography}{99}

\par\bibitem{AG}    J. T. Arnold and R. Gilmer, On the contents of polynomials, Proc. Amer. Math. Soc. {\bf 24} (1970), 556--562.
\par\bibitem{BM}    C. Bakkari and N. Mahdou, Gaussian polynomials and content ideal in pullbacks, Comm. Algebra {\bf 34} (8) (2006), 2727--2732.
\par\bibitem{BG1}    S. Bazzoni and S. Glaz, Pr\"ufer rings in ``Multiplicative Ideal Theory in Commutative Algebra" Eds:
J. W. Brewer, S. Glaz, W. J. Heinzer, and B. Olberding, Springer, pp. 263--277, 2006.
\par\bibitem{BG2}    S. Bazzoni and S. Glaz, Gaussian properties of total rings of quotients, J. Algebra {\bf 310} (2007), 180--193.
\par\bibitem{BS1}     M. Boisen and P. Sheldon, A note on pre-arithmetical rings, Acta Math. Acad. Sci. Hungar. {\bf 28} (1976), 257--259.
\par\bibitem{BS2}     M. Boisen and P. Sheldon, Pre-Pr\"ufer rings, Pac. J. Math., {\bf 58} (1975), 331--344.
\par\bibitem{BS}    H. S. Butts and W. Smith, Pr\"ufer rings, Math. Z. {\bf 95} (1967), 196--211.
\par\bibitem{GH}    R. Gilmer and J. Huckaba,  $\Delta $-rings, J. Algebra {\bf 28} (1974), 414--432.
\par\bibitem{G1}    S. Glaz, Commutative coherent rings, Lecture Notes in Mathematics, 1371, Springer-Verlag, Berlin, 1989.
\par\bibitem{G2}    S. Glaz, The weak dimension of Gaussian rings, Proc. Amer. Math. Soc. {\bf 133} (9) (2005), 2507--2513.
\par\bibitem{G3}    S. Glaz, Pr\"ufer conditions in rings with zero-divisors, CRC Press Series of Lectures in
                 Pure Appl. Math. {\bf 241} (2005), 272--282.
\par\bibitem{GV}    S. Glaz and W. Vasconcelos, The content of Gaussian polynomials, J. Algebra {\bf 202} (1998), 1--9.
\par\bibitem{Gr}     M. Griffin, Pr\"ufer rings with zero-divisors, J. Reine Angew Math. {\bf 239/240} (1970), 55--67.
\par\bibitem{HH1}   W. Heinzer and C. Huneke, Gaussian polynomials and content ideals, Proc. Amer. Math. Soc. {\bf 125} (1997), 739--745.
\par\bibitem{HH2}   W. Heinzer and C. Huneke, The Dedekind-Mertens lemma and the content of polynomials, Proc. Amer. Math. Soc. {\bf 126} (1998), 1305--1309.
\par\bibitem{H}     J. A. Huckaba, Commutative rings with zero divisors, Marcel Dekker, New York, 1988.
\par\bibitem{J}     C. U. Jensen, Arithmetical rings, Acta Math. Hungr. {\bf 17} (1966), 115--123.
\par\bibitem{KM}   S. Kabbaj and N. Mahdou, Trivial extensions defined by coherent-like conditions, Comm. Algebra {\bf 32} (10) (2004), 3937--3953.
\par\bibitem{Kr}    W. Krull, Beitrage zur arithmetik kommutativer integritatsbereiche, Math. Z. {\bf 41} (1936), 545--577.
\par\bibitem{LR}    K. A. Loper and M. Roitman, The content of a Gaussian polynomial is invertible, Proc. Amer. Math. Soc. {\bf 133} (5) (2004), 1267--1271.
\par\bibitem{Lu}   T. G. Lucas, Gaussian polynomials and invertibility, Proc. Amer. Math. Soc. {\bf 133} (7) (2005), 1881--1886.
 \par\bibitem{P}     H. Pr\"ufer, Untersuchungen uber teilbarkeitseigenschaften in korpern, J. Reine Angew. Math. {\bf 168} (1932), 1--36
\par\bibitem{Ro}    J. J. Rotman, An introduction to homological algebra, Academic Press, New York, 1979.
\par\bibitem{Ru}    D. E. Rush, The Dedekind-Mertens lemma and the contents of polynomials, Proc. Amer. Math. Soc. {\bf 128} (2000), 2879--2884.
\par\bibitem{T}     H. Tsang, Gauss's Lemma, Ph.D. thesis, University of Chicago, Chicago, 1965.
\end{thebibliography}
\end{document}